\documentclass[11pt,a4paper]{report}

\usepackage[margin=2.6cm]{geometry}
\usepackage{amsmath,amssymb,amsthm}
\usepackage{enumitem}
\usepackage{xcolor}
\usepackage[colorlinks=true,linkcolor=blue!60!black,citecolor=blue!60!black,urlcolor=blue!60!black]{hyperref}

\theoremstyle{plain}
\newtheorem{theorem}{Theorem}[chapter]
\newtheorem{lemma}[theorem]{Lemma}

\newtheorem{proposition}[theorem]{Proposition}

\theoremstyle{definition}
\newtheorem{definition}[theorem]{Definition}

\theoremstyle{remark}
\newtheorem{remark}[theorem]{Remark}

\numberwithin{equation}{chapter}

\newcommand{\E}{\mathbb{E}}
\newcommand{\Prob}{\mathbb{P}}
\newcommand{\Var}{\mathrm{Var}}
\newcommand{\Z}{\mathbb{Z}}

\newcommand{\N}{\mathbb{N}}
\newcommand{\Reff}{R^{\mathrm{eff}}}

\newcommand{\iio}{\ \text{i.o.}}
\newcommand{\as}{\ \text{a.s.}}
\newcommand{\ind}[1]{\mathbf{1}_{\{#1\}}}

\title{\bfseries Random Walk on age-preferential random growth processes}
\author{Jia-shan Tang,\  Feng Wang,\  Xian-Yuan Wu}
\date{}

\begin{document}
\maketitle

\begin{abstract}
This paper studies a class of \textbf{exponentially age-preferential random growth graphs} $G(m,\beta)$: at each time $n$ a new vertex $v_n$ is added, and $m$ edges are drawn to existing vertices independently, with the probability that a target $v_k$ is chosen proportional to $e^{\beta k}$ ($\beta>0$), so that the graph's growth strongly favors ``young'' vertices. We systematically investigate two fundamental probabilistic problems of this model and supply the proof details: (1)~\textbf{Recurrence}: when $\beta>0$, the simple random walk on $G$ is recurrent with probability $1$ (Theorem~\ref{thm:rec}); (2)~\textbf{Infinite collisions}: for every $\beta>0$, two independent simple random walks on $G$ meet infinitely often with probability $1$ (Theorem~\ref{thm:icp}). The core of the proof is the following: by the Borel--Cantelli lemma combined with the \textbf{memoryless property} of exponential weights, $G$ is almost surely ``sausage-like'' --- there exist infinitely many bottlenecks consisting of exactly $m$ parallel edges; the cascade of bottlenecks makes the effective resistance diverge, yielding recurrence; the renewal process $\{T_k\}$ constructed at the bottleneck times has i.i.d.\ stationary increments with finite moments of every order for the first-passage times, and combined with a moderate-deviation estimate we verify Peres' Green-function criterion for infinite collisions.

\vspace{4pt}
\end{abstract}

\tableofcontents

\chapter{Introduction}\label{ch:intro}
\section{Background}\label{sec:bac}
Random walks on random graphs form an important research topic at the intersection of probability and combinatorics. Given a locally finite connected graph $G=(V,E)$, the simple random walk on $G$ is the Markov chain which starts at a given vertex and, at each step, moves with equal probability to one of the neighbours of the current vertex. The most fundamental qualitative problem is the recurrence of the walk: whether the walk returns to its starting vertex with probability $1$. The classical P\'olya theorem shows that the simple random walk on the integer lattice $\Z^d$ is recurrent when $d\leq 2$ and transient when $d\geq 3$. The electric-network theory developed by Doyle and Snell \cite{DS} and by Lyons and Peres \cite{LP} gives a unified criterion: $G$ is recurrent if and only if the effective resistance from any vertex to infinity, $\Reff(v_0,\infty)=\infty$. This criterion transforms the recurrence of a random walk into the computation or estimation of effective resistance, and is one of the main tools of this paper. For randomly generated graphs, Benjamini and Schramm \cite{BS} studied recurrence of distributional limits of finite planar graphs, and Gurel-Gurevich and Nachmias \cite{GGN} further proved that distributional limits of recurrent graphs remain recurrent.

On top of recurrence, a finer question is the \textbf{infinite collision property} of independent random walks. Let $G$ be an infinite connected graph, and let $ X = \{X_n\} $ and $ X' = \{X'_n\} $ be two independent simple random walks on $G$, we say that $G$ has the \textit{infinite collision property} if $ | \{ n : X_n = X'_n \} | = \infty $ (i.e. occupy the same vertex at the same time) almost surely. Otherwise, we say $G$ has the finite collision property. Not every recurrent graph has the infinite-collision property:  Krishnapur and Peres \cite{KP} constructed recurrent graphs---a comb lattice and its variants---on which two independent walks collide only finitely often, so recurrence does not force ICP. Barlow, Peres and Sousi \cite{BPS} gave a sufficient condition through the Green function: on a recurrent graph, if some exhaustion by finite subgraphs keeps $\sup_x g(x,x)/g(\rho,\rho)$ bounded, then ICP holds; they also obtained ICP for $\Z^d$ with $d\le4$, for recurrent trees, and for uniform spanning forests. Later work tested the criterion on specific geometries. Chen and Chen \cite{CC1} proved ICP for two walks on the open cluster of $\Z^2$, and in \cite{CC2} gave sufficient conditions on wedge combs, later extended to general wedges by Chen \cite{C}. Hutchcroft and Peres \cite{HP} proved that every recurrent reversible random rooted graph has ICP, covering unimodular random graphs and spanning trees; reversibility is essential. Halberstam and Hutchcroft \cite{HH} opened the dynamic setting with collisions in the dynamic random conductance model on $\Z^2$. Recent contributions diversify in three directions. Quantitatively, Watanabe \cite{Wa} estimated collision times for two walks on the uniform spanning tree of $\Z^3$, giving a new proof of ICP that bypasses Green functions; Croydon, Shiraishi and Watanabe \cite{CSW} moved to random walk traces, studying collisions on the trace of a four-dimensional walk. Structurally, De Ambroggio, Nitzschner and Scali \cite{DNS} located phase transitions for comb graphs with planar bases. Probabilistically, Astoquillca \cite{As1} used the duality between the voter model and coalescing random walks to link ICP with stationary measures of interacting particle systems, and in \cite{As2} proved that discrete- and continuous-time ICP agree on bounded-degree graphs.

In this paper we consider the \textbf{age-preferential random growth process}: unlike the classical Barab\'asi--Albert degree-preferential model, here the probability with which a new edge selects its target depends only on the \textbf{label} of the target (i.e., the complement of its ``age''), and increases with the label --- the younger the vertex, the more likely it is to be chosen. In the earlier work \cite{TWW1,TWW2}, Tang, Wang and Wu studied the bounded-preference case (new vertices may connect only to the most recent finite number of vertices) and random recursive trees generated by Bernoulli sequences; this paper focuses on the exponential case, fixing $m\geq 1$, with each time step adding one new vertex and $m$ independent edges to existing vertices, the parent being chosen proportionally to the exponential weight $e^{\beta k}$ ($\beta>0$), i.e., the $k$-th-oldest vertex is selected with probability proportional to $e^{\beta k}$. By the memoryless property of the exponential distribution, this model enjoys a translational self-similarity, which is the root of the renewal theory we are able to build.
\section{Main results}\label{sec:res}
Let $G(m,\beta)$ be the limit graph of the above exponentially age-preferential random growth process. The main results proven in this paper are as follows.

\begin{theorem}[Recurrence]\label{thm:rec}
For every $m\geq 1$ and $\beta>0$, the limit graph $G(m,\beta)$ is recurrent $\as$.
\end{theorem}

\begin{theorem}[Infinite collisions]\label{thm:icp}
For every $m\geq 1$ and $\beta>0$, $G(m,\beta)$ has the infinite-collision property $\as$; moreover the number of times that two independent simple random walks, at the same time step, traverse the same edge is also infinite $\as$.
\end{theorem}

The proofs of the two theorems above rely on a common geometric structure of the graph: almost surely, $G$ has a \textbf{sausage-like} structure --- there exist infinitely many times such that any path through the graph must go through a ``bottleneck'' consisting of exactly $m$ parallel edges. Each bottleneck contributes at least $1/m$ resistance in series, making the effective resistance diverge and yielding recurrence (Theorem~\ref{thm:rec}); the bottleneck also gives an estimate showing that the Green function at the root grows linearly in the graph size, which verifies Peres' criterion and yields the infinite-collision property (Theorem~\ref{thm:icp}). In the proof of Theorem~\ref{thm:icp} we also establish a key renewal structure: by the memoryless property of the exponential distribution, the ``bottleneck times'' form a renewal process whose increments are i.i.d.\ with finite moments of every order; this is the technical core of our proof (see Theorem 

\section{Outline}\label{sec:rest}
The rest of the paper is organized as follows. Chapter~2 is preliminaries: \S\ref{sec:model} gives a rigorous definition of the model and the notation used throughout; \S\ref{sec:prelim} collects the basic facts of electric-network theory and probability used throughout; \S\ref{sec:lem} states and briefly proves the literature results that will be quoted in later chapters. Chapter~3 studies recurrence: \S\ref{sec:saus} establishes the sausage-like structure of the graph (tail estimates, time scales, and a lower bound on the probability of bottleneck events), and \S\ref{sec:proofrec} completes the proof of Theorem~\ref{thm:rec}. Chapter~4 studies infinite collisions: \S\ref{sec:renewal} defines the renewal times and proves the renewal theorem (Theorem 
), \S\ref{sec:moments} provides moment estimates of the renewal lengths, and \S\ref{sec:proof-icp} completes the proof of Theorem~\ref{thm:icp} via moderate deviations and Peres' criterion. Chapter~5 summarizes the paper and gives some outlooks including the power-law preference case.

\chapter{Preliminaries}\label{ch:prem}

\section{Model definition and notation}\label{sec:model}

\begin{definition}[Exponentially age-preferential random growth process]\label{def:model}
Fix an integer $m\ge 1$ and a parameter $\beta>0$. Let $\{G_n=(V_n,E_n):n\ge 0\}$ be a random graph process such that $G_0=(V_0,E_0)$, $V_0=\{v_0\}$, $E_0=\varnothing$. The rule for going from $G_{n-1}$ to $G_n$ ($n\ge 1$) is:

\begin{enumerate}
  \item[(1)] $V_n = V_{n-1}\cup\{v_n\}$, i.e., a new vertex $v_n$ is added;
  \item[(2)] $E_n = E_{n-1}\cup\{\langle v_n, w_{n,i}\rangle:1\le i\le m\}$, where the random variables $w_{n,1},\dots,w_{n,m}\in V_{n-1}$ are i.i.d., and
  \begin{equation}\label{eq:model}
    \Prob(w_{n,i}=v_k)=\frac{e^{\beta k}}{Z_n},
    \qquad k=0,1,\dots,n-1,
  \end{equation}
  where the normalizing constant is
  \begin{equation}\label{eq:Zn}
    Z_n=\sum_{k=0}^{n-1}e^{\beta k}=\frac{e^{\beta n}-1}{e^{\beta}-1}.
  \end{equation}
\end{enumerate}
Define the limit graph
\[
  G=\bigcup_{n=0}^{\infty}G_n=:\,G(m,\beta),\qquad
  V=\bigcup_{n=0}^{\infty}V_n=\{v_0,v_1,v_2,\dots\}.
\]

\end{definition}

\begin{remark}
(i)~ $G$ is in general a multigraph: there may be multiple parallel edges between $v_n$ and the same parent. Throughout this paper ``graph'' allows parallel edges but no self-loops; the notions of random walk and spanning tree are defined for multigraphs in the usual way, with transition probabilities counted by the number of parallel edges.

(ii)~When the vertex $v_n$ joins the model, the age of $v_k$ is $n-k$. In the distribution \eqref{eq:model}, the weight $e^{\beta k}$ increases as the vertex age decreases (the smaller the index $k$, the older the vertex---so the weight is smallest for the oldest vertices and largest for the most recent ones). Equivalently, new vertices preferentially attach to the most recently added vertices, which is the meaning of the term ``age preference'' in our terminology.
\end{remark}

\begin{definition}[Random walk]\label{def:srw}
Let $H=(V_H,E_H)$ be a locally finite connected graph. The simple random walk $\{X_n:n\geq 0\}$ on $H$ is the discrete-time Markov chain on the state space $V_H$ with transition probabilities
\[
p_{u,v}=\begin{cases}\dfrac{\#\{e\in E_H:\,e=\{u,v\}\}}{\deg_H(u)},& v\sim u,\\[2mm] 0,&\text{otherwise},\end{cases}
\]
where $\deg_H(u)$ is the degree of $u$ in $H$.

If $\Prob_x(X_n=x \text{ infinitely often})=1$ for every $x$, then $H$ (the walk on $H$) is said to be \textbf{recurrent}. Naturally, if $\{Y_t: t\in\mathbb{Z}_+\}$ is a continuous-time random walk on $H$, its transition density is defined by
\[
q_t(x,y)=\frac{\mathbb{P}^x(Y_t=y)}{\deg_H(u)}.
\]
\end{definition}

\begin{definition}[Green kernel]\label{def:Green}
For a locally finite connected graph $G$, let $Y_t$ be the random walk on $G$ and $q_t$ its transition density. The \textbf{Green kernel} is defined as
\[
g_G(x,y)=\sum_{t=0}^{\infty} q_t(x,y).
\]
We generally abbreviate $g_G(x,y)$ as $g(x,y)$. For a finite subset $H$ of $G$, the \textbf{killed Green kernel} is defined as
\[
g_H(x,y)=\sum_{t=0}^{\tau_{H^c}} q_t(x,y).
\]
where $\tau_{H^c}$ is the first time the walk leaves $H$.
\end{definition}

\begin{definition}[Finite and infinite collisions]\label{def:collision}
Let $X,Y$ be two independent simple random walks on a graph $G$. Let $\Prob_{u,v}$ denote the distribution of the process $\{(X_n,Y_n):n\geq 0\}$ with $X_0=u,Y_0=v$. Define
\[
Z:=\sum_{n=0}^{\infty}\ind{X_n=Y_n}
\]
as the total number of meetings. If $\Prob_{u,u}(Z=\infty)=1$ for all $u\in V$, we say that $G$ has the \textbf{infinite-collision property} (ICP); if $\Prob_{u,u}(Z<\infty)=1$, then $G$ has the \textbf{finite-collision property}.
\end{definition}

The other commonly used notation throughout the rest of the paper is as follows. For a realization $G=G(m,\beta)$, let $d_G(\cdot,\cdot)$ be the graph distance; for $r\geq 1$, $B_r:=\{u\in V:d_G(v_0,u)\leq r\}$. $\Reff_H(x\leftrightarrow y)$ and $\Reff_H(x\leftrightarrow H^c)$ denote respectively the effective resistance in the network $H$ between $x$ and $y$, and between $x$ and the complement of $H$ (which we identify with a single vertex); we abbreviate $\Reff_G(x\leftrightarrow y)$ as $\Reff(x\leftrightarrow y)$. For a sequence of events $\{A_n\}$, $A_n\iio$ (``infinitely often'') means that $A_n$ occurs infinitely many times. Fix $\alpha,\delta>0$, let $\lfloor x\rfloor$ denote the floor function, and abbreviate $V_{\lfloor n-\alpha\ln n\rfloor}$ as ``the set of all vertices whose indices do not exceed $n-\alpha\ln n$''.

\section{Preliminaries}\label{sec:prelim}

We first collect the basic facts of electric-network theory; for proofs see Doyle--Snell \cite{DS} and Lyons--Peres \cite{LP}.

\begin{proposition}[Basic facts of electric networks]\label{prop:electric}
Let $G$ be a locally finite, infinite, connected network (each edge having resistance $1$), and $x,y\in H \subset G$. Then:

(i)~ (Green-function / effective-resistance identity) For $H$, $g_H(x,x)=\Reff_H(x\leftrightarrow H^c)$; similarly, $g(x,x)=\Reff(x\leftrightarrow\infty)$.

(ii)~ (Recurrence criterion) $G$ is recurrent iff $\Reff(x\leftrightarrow\infty)=\infty$ for some (and hence every) $x$; equivalently, $g(x,x)$ diverges along an exhausting sequence:
\[
 G \text{ is recurrent}\iff g(x,x)=\lim_{r\to\infty}g_{B_r}(x,x)=\Reff(x\leftrightarrow\infty)=+\infty .
\]

(iii)~ (Rayleigh monotonicity) Adding (resp.~deleting) edges of the network can only decrease (resp.~increase) the effective resistance; short-circuiting several vertices can only decrease the effective resistance.

(iv)~ (Energy principle / Thomson's principle) $\Reff_H(x\leftrightarrow y)$ equals the minimum energy of a unit flow from $x$ to $y$; in particular, if edge cuts $C_1,\dots,C_J$ separate $x$ from the target set $S$ \emph{in series} (i.e., every $x$--$S$ path must traverse each cut in turn), then
\begin{equation}\label{eq:cutbound}
\Reff_H(x\leftrightarrow S)\;\geq\;\sum_{j=1}^{J}\frac{1}{|C_j|},
\end{equation}
where $|C_j|$ is the number of edges in the cut $C_j$.

(v)~ For two distinct points $x\neq y$ in a finite network, $\Reff_H(x\leftrightarrow y)\geq 1/|E_H|$.
\end{proposition}

Next we recall facts related to the Borel--Cantelli lemma and classical limit theorems of probability; for proofs see Durrett \cite{Du} and Dembo--Zeitouni \cite{DZ}.

\begin{proposition}[Borel--Cantelli lemma]\label{prop:bc}
Let $\{A_n\}$ be a sequence of events.

(i)~ If $\sum_n \Prob(A_n)<\infty$, then $\Prob(A_n\iio)=0$.

(ii)~ If $\{A_n\}$ are independent and $\sum_n\Prob(A_n)=\infty$, then $\Prob(A_n\iio)=1$.
\end{proposition}

\begin{proposition}[Limit theorems for sums of i.i.d.\ random variables]\label{prop:limit}
Let $\{X_n:n\geq 1\}$ be i.i.d., $S_n=\sum_{k=1}^n X_k$. When the corresponding moment conditions hold:

(i)~ (Strong law of large numbers) If $\E|X_1|<\infty$, $\mu=\E X_1$, then $S_n/n\to\mu\as$.

(ii)~ (Central limit theorem) If $\E X_1=\mu$, $0<\sigma^2=\Var(X_1)<\infty$, then $(S_n-n\mu)/(\sigma\sqrt n)\Rightarrow N(0,1)$.

(iii)~ (Moderate deviations, MDP) If $\E X_1=\mu$, $\Var(X_1)=\sigma^2\in(0,\infty)$, and $\{x_n\}$ satisfies $x_n\to\infty$, $x_n=o(n)$, then
\begin{equation}\label{eq:mdp}
\Prob\bigl(S_n-n\mu\geq x_n\bigr)\leq \exp\Bigl(-\frac{x_n^2}{2n\sigma^2}\,(1+o(1))\Bigr).
\end{equation}
In particular, taking $x_n=n^{3/4}$ gives $\Prob(S_n-n\mu\geq n^{3/4})\leq \exp(-c\,n^{1/2})$ (for $n$ sufficiently large, $c=c(\sigma)>0$), and this bound is summable in $n$.
\end{proposition}

\section{Peres' criterion and a short proof}\label{sec:lem}
\begin{lemma}[Peres' criterion, Barlow--Peres--Sousi \cite{BPS}]\label{lem:peres}
Let $G$ be a recurrent graph and $o\in V$ a fixed vertex. Let $(B_r)_{r}$ be an increasing sequence of subsets of $V$ such that $B_r\neq G$ for every $r$ and $\bigcup_r B_r=G$. If there exists a constant $C<\infty$ such that for every $r$,
\[
g_{B_r}(x,x)\;\le\;C\, g_{B_r}(o,o),\qquad \forall\, x\in B_r,
\]
then $G$ has the infinite-collision property. Moreover, the total number of times that two independent random walks, at the same time step, traverse the same edge, is almost surely infinite.
\end{lemma}

\begin{proof}
Let $B:=B_r$, and let $X^B$ and $Y^B$ be two independent random walks started from $o$ and killed when leaving $B$, with transition density $q_t^B$. Define
\[
\widetilde{Z}^B=\sum_{t=0}^{\infty}
\mathbf{1}\bigl(X_t^B=Y_t^B,\; X_{t+1}=Y_{t+1}\bigr)
\]
as the total number of times the two walks, at the same time step, traverse the same edge. The idea of the proof is to apply the second-moment method to $\widetilde{Z}^B$.

\textbf{First moment.} Using the independence of the two walks and summing directly (the degree normalizing factors cancel when summing over $x\in B$ and $y\sim x$), we obtain
\[
\mathbb{E}_{o,o}\bigl[\widetilde{Z}^B\bigr]
=\sum_{t=0}^{\infty} q_{2t}^{B}(o,o)
\;\in\;\Bigl[\tfrac12\,g_B(o,o),\; g_B(o,o)\Bigr],
\]
where we used the fact that $G$ is recurrent and $B_r\neq G$, so that $g_{B_r}(o,o)<\infty$.

\textbf{Second moment.} Decompose $\mathbb{E}[(\widetilde{Z}^B)^2]$ according to ``first cross edge $(x,y)$ at time $t$, then cross edge $(z,w)$ at time $s>t$''. Summing over $z,w$ and controlling the remaining sum by the Green kernel gives
\[
\mathbb{E}_{o,o}\bigl[(\widetilde{Z}^B)^2\bigr]
\;\le\; g_B(o,o)+2\,g_B(o,o)\max_{y\in B} g^B(y,y).
\]

\textbf{Second-moment method and passage to the limit.} By the Paley--Zygmund (second-moment) inequality, together with the two estimates above and the hypothesis $\max_{y\in B_r} g_{B_r}(y,y)\le C\, g_{B_r}(o,o)$, we obtain
\[
\mathbb{P}_{o,o}\Bigl(\widetilde{Z}^{B_r}>\tfrac12\,\mathbb{E}_{o,o}[\widetilde{Z}^{B_r}]\Bigr)
\;\ge\;\frac{g_{B_r}(o,o)}{16\bigl(1+2C\,g_{B_r}(o,o)\bigr)}.
\]
Since $g_{B_r}(o,o)\ge d(o)^{-1}$, there is a constant $c>0$, independent of $r$, such that
\[
\mathbb{P}_{o,o}\Bigl(\widetilde{Z}^{B_r}>\tfrac14\, g_{B_r}(o,o)\Bigr)\;\ge\; c,
\qquad \forall r.
\]
Letting $r\to\infty$, we have $\widetilde{Z}^{B_r}\uparrow\widetilde{Z}$ (where $\widetilde{Z}$ is the total number of edges the two walks traverse in common on $G$), and by monotone convergence $\mathbb{P}_{o,o}(\widetilde{Z}=\infty)\ge c$. Since the total number of vertex-collisions $Z\ge\widetilde{Z}$, we have $\mathbb{P}_{o,u}(Z=\infty)\ge c$; applying the 0--1 law \cite{BPS} yields $\mathbb{P}_{o,o}(Z=\infty)=1$, i.e., $G$ has the infinite-collision property.
\end{proof}

\chapter{Recurrence}\label{ch:rec}

This chapter proves Theorem~\ref{thm:rec}. The strategy is: first, prove that $G$ is almost surely ``sausage-like'', i.e., contains infinitely many bottlenecks consisting of $m$ parallel edges; then use the energy principle (Proposition~\ref{prop:electric}(iv)) to accumulate infinitely many bottlenecks into a divergent effective resistance; finally appeal to the recurrence criterion (Proposition~\ref{prop:electric}(ii)).

\section{Sausage-like structure}\label{sec:saus}
\begin{lemma}[Tail estimate]\label{lem:tail}
Let $\alpha>1/\beta$. For $n\geq 1$, define
\[
D_n:=\Bigl\{\exists\, i\in\{1,\dots,m\}:\ w_{n,i}\in V_{\lfloor n-\alpha\ln n\rfloor}\Bigr\},
\]
i.e., some new edge connects to an old vertex of index $\leq n-\alpha\ln n$. Then
\[
\Prob(D_n)\leq m e^{2\beta}\,n^{-\alpha\beta},
\qquad\text{hence}\qquad
\sum_{n=1}^{\infty}\Prob(D_n)<\infty.
\]
In particular, by the Borel--Cantelli lemma (Proposition~\ref{prop:bc}(i)), almost surely, for all sufficiently large $n$ and all $1\leq i\leq m$,
\begin{equation}\label{eq:tailevent}
\text{the index of $w_{n,i}$}\;>\;n-\alpha\ln n .
\end{equation}
\end{lemma}

\begin{proof}
By the union bound,
\[\Prob(D_n)
 \le m\,\Prob(w_{n,1}\in V_{n-\alpha\ln n})
\]
and by \eqref{eq:model}
\[
\Prob(w_{n,1}\in V_{n-\alpha\ln n})\leq \,\frac{Z_{\lfloor n-\alpha\ln n+1\rfloor}}{Z_n}
=\,\frac{e^{\beta\lfloor n-\alpha\ln n+1\rfloor}-1}{e^{\beta n}-1}
\leq \,e^{\beta(n-\alpha\ln n+1)-\beta (n-1)}
= \,e^{2\beta}\,n^{-\alpha\beta}.
\]
Since $\alpha\beta>1$, the right-hand side is summable in $n$.

Lemma~\ref{lem:tail} shows that new edges almost surely connect only to the $\alpha\ln n$ most recent vertices; i.e., the graph's growth is local in time. To construct bottlenecks we now pick a sequence of logarithmically-spaced times.

\begin{lemma}[Logarithmic time scale]\label{lem:scale}
Let $\delta>1$, $n_0$ sufficiently large, and define recursively
\begin{equation}\label{eq:nk}
n_{k+1}=n_k+\lceil \delta\ln n_k\rceil,\qquad k\geq 0 .
\end{equation}
Then there exist a constant $M_0<\infty$ depending only on $(\delta,n_0)$ such that for all $k\geq 1$:

(i)~ $n_k\geq n_0+\delta k\ln n_0$.

(ii)~ $\dfrac{\ln n_k}{n_k}\leq \dfrac1k$, and hence $n_k\leq n_0 e^{\delta}\,k^{\delta}\leq M_0 k^{\delta}$.
\end{lemma}

\begin{proof}
(i)~ From the recursion \eqref{eq:nk} and $n_k\ge n_0$ we have
\[
n_{k+1}-n_k=\lceil\delta\ln n_k\rceil\ge \delta\ln n_k\ge \delta\ln n_0.
\]
Summing for $k$ from $0$ to $k-1$ gives
\[
n_k\ge n_0+\delta k\ln n_0.
\]
In particular $n_k\to\infty$ (uniformly, so long as $n_0>1$).

(ii)~ We prove (ii) by induction. Define
\[
f(x)=\frac{\ln x}{x},\qquad x>0.
\]
Direct differentiation gives
\[
f'(x)=\frac{1-\ln x}{x^2}.
\]
Hence $f$ is strictly decreasing for $x>e^{3/2}$.

For $k=1$: since $n_1=n_0+\lceil\delta\ln n_0\rceil>n_0$, we have
\[
\frac{\ln n_1}{n_1}=f(n_1)<f(n_0)=\frac{\ln n_0}{n_0}.
\]
Choosing $n_0$ large enough so that $\dfrac{\ln n_0}{n_0}\le 1$ (e.g.\ $n_0\ge 3$) gives
\[
\frac{\ln n_1}{n_1}\le 1.
\]

Assume inductively that for some $k\ge 1$
\begin{equation}\label{eq:IH}
\frac{\ln n_k}{n_k}\le \frac1k,
\qquad\text{i.e.,}\qquad
n_k\ge k\ln n_k.
\end{equation}

We need to prove
\[
\frac{\ln n_{k+1}}{n_{k+1}}\le \frac1{k+1},
\]
which is equivalent to
\begin{equation}\label{eq:goal}
n_{k+1}\ge (k+1)\ln n_{k+1}.
\end{equation}

Define
\[
g_k(x):=x-(k+1)\ln x,\qquad x>0.
\]
Then \eqref{eq:goal} says $g_k(n_{k+1})\ge 0$. We now estimate $g_k(n_{k+1})$ in terms of $n_k$.

From the recursion \eqref{eq:nk},
\[
\ln n_{k+1}=\ln\Bigl(n_k\Bigl(1+\frac{\lceil\delta\ln n_k\rceil}{n_k}\Bigr)\Bigr)=\ln n_k+\ln\!\left(1+\frac{\lceil\delta\ln n_k\rceil}{n_k}\right).
\]
Therefore
\begin{align}
g_k(n_{k+1})
&=n_k+\lceil\delta\ln n_k\rceil-(k+1)\ln n_k-(k+1)\ln\!\left(1+\frac{\lceil\delta\ln n_k\rceil}{n_k}\right)\notag\\[4pt]
&=\bigl(n_k-k\ln n_k\bigr)+\bigl(\lceil\delta\ln n_k\rceil-\ln n_k\bigr)
-(k+1)\ln\!\left(1+\frac{\lceil\delta\ln n_k\rceil}{n_k}\right).\label{eq:gk}
\end{align}

By the induction hypothesis \eqref{eq:IH}, the first term $n_k-k\ln n_k\ge 0$. Now we control the next two terms. Set
\[
L:=\ln n_k,\qquad N:=n_k.
\]
By \eqref{eq:IH} we have $N\ge kL$, and
\[
\frac{\lceil\delta L\rceil}{N}
\le \frac{\delta L+1}{N}
=\delta\,\frac{L}{N}+\frac1N
\le \frac{\delta}{k}+\frac1N
\le \frac{\delta}{k}+\frac1{n_0}.
\]
When $k$ or $n_0$ is large enough, the right-hand side is $<1$, so by $\ln(1+x)\le x$ ($x>-1$),
\[
\ln\!\left(1+\frac{\lceil\delta L\rceil}{N}\right)
\le \frac{\lceil\delta L\rceil}{N}
\le \frac{\delta L+1}{N}.
\]
Plugging into \eqref{eq:gk} and using $N\ge kL$,
\begin{align*}
g_k(n_{k+1})
&\ge \lceil\delta L\rceil-L-(k+1)\frac{\delta L+1}{N}\\[4pt]
&\ge \delta L-L-(k+1)\frac{\delta L+1}{kL}\\[4pt]
&=(\delta-1)L-\frac{k+1}{k}\left(\delta+\frac1L\right).
\end{align*}

Since $\delta>1$ and by (i) we have $L=\ln n_k\to\infty$ (uniformly for $k$ as $n_0$ is large), there exists a threshold $L_0(\delta)$ depending only on $\delta$ such that whenever $L\ge L_0(\delta)$,
\[
(\delta-1)L\ge 2\left(\delta+1\right)
\ge \frac{k+1}{k}\left(\delta+\frac1L\right)
\qquad(\forall\,k\ge 1).
\]
And $L=\ln n_k\ge L_0(\delta)$ holds so long as $n_0$ is large enough (since $n_k\ge n_0$). Therefore
\[
g_k(n_{k+1})\ge 0,
\]
which is \eqref{eq:goal}. The induction is complete, and we have for all $k\ge 1$
\[
\frac{\ln n_k}{n_k}\le \frac1k.
\]

Taking logarithms in the recursion \eqref{eq:nk},
\begin{equation}\label{eq:lnnk}
\ln n_{k+1}-\ln n_k
=\ln\!\left(1+\frac{\lceil\delta\ln n_k\rceil}{n_k}\right)
\le \frac{\lceil\delta\ln n_k\rceil}{n_k}
\le \delta\,\frac{\ln n_k}{n_k}+\frac1{n_k}.
\end{equation}
Summing for $i$ from $1$ to $k$, and using (ii) together with $n_i\ge n_0$, we get
\[
\ln n_k
\le \ln n_0+\delta\sum_{i=1}^{k}\frac1i+\sum_{i=1}^{k}\frac1{n_0}.
\]
Since $\sum_{i=1}^{k}\frac1i\le \ln k+1$, we obtain
\[
\ln n_k\le \ln n_0+\delta\ln k+\delta+\frac{k}{n_0}=\ln n_0+\delta\ln k+O(1).
\]
With $n_0\ge 1$ we conclude
\[
n_k\le M_0\,k^{\delta},
\]
where $M_0=C(\delta)n_0$ is some constant depending only on $(\delta,n_0)$.
\end{proof}

We now introduce the central event of this chapter. Fix $1/\beta<\alpha<\delta$ and the sequence of times \eqref{eq:nk}, and define the \textbf{sausage event}
\begin{equation}\label{eq:Ak}
A_k:=\Bigl\{\begin{array}{l}
\text{(a)}\ \forall\,1\leq i\leq m:\ w_{n_k+1,i}=v_{n_k};\\[1mm]
\text{(b)}\ \forall\, n_k+1<n\leq n_{k+1},\ \forall\,1\leq i\leq m:\ w_{n,i}\in\{v_{n_k+1},\dots,v_{n-1}\}.
\end{array}\Bigr.
\end{equation}
On $A_k$, all $m$ edges from $v_{n_k+1}$ go to $v_{n_k}$ (forming a bottleneck), and the whole subsequent block of new vertices $\{v_{n_k+1},\dots,v_{n_{k+1}}\}$ only has internal edges --- a ``sausage'' hanging from the bottleneck.

\begin{lemma}[Lower bound and independence of sausage events]\label{lem:sausage}
With the notation above,
\[
\Prob(A_k)\;\geq\;c(\beta):=(1-e^{-\beta})^{m}\Bigl(\prod_{s=1}^{\infty}\bigl(1-e^{-\beta s}\bigr)\Bigr)^{m}\;>\;0,
\]
where the infinite product converges to a positive number because $\sum_s e^{-\beta s}<\infty$. Moreover $\{A_k:k\geq 0\}$ are mutually independent, and consequently $\Prob(A_k\iio)=1$.
\end{lemma}

\begin{proof}
Since $Z_{n_k+1}=(e^{\beta(n_k+1)}-1)/(e^{\beta}-1)$,
\[
\frac{e^{\beta n_k}}{Z_{n_k+1}}
=\frac{e^{\beta n_k}(e^{\beta}-1)}{e^{\beta(n_k+1)}-1}
\geq \frac{e^{\beta}-1}{e^{\beta}}=1-e^{-\beta}.
\]

For $n_k+1<n\leq n_{k+1}$,
\[
\Prob\bigl(w_{n,i}\in\{v_{n_k+1},\dots,v_{n-1}\}\bigr)
=1-\frac{Z_{n_k+1}}{Z_n}
=\frac{e^{\beta n}-e^{\beta(n_k+1)}}{e^{\beta n}-1}
\geq 1-e^{-\beta(n-n_k-1)} .
\]
By independence of the choices at distinct $(n,i)$,
\[
\Prob(A_k)\geq (1-e^{-\beta})^{m}\prod_{n=n_k+2}^{n_{k+1}}\bigl(1-e^{-\beta(n-n_k-1)}\bigr)^{m}
\geq (1-e^{-\beta})^{m}\Bigl(\prod_{s=1}^{\infty}(1-e^{-\beta s})\Bigr)^{m}.
\]
Independence: $A_k$ depends only on $\{w_{n,i}:n_k<n\leq n_{k+1}\}$; the intervals $(n_k,n_{k+1}]$ are pairwise disjoint and the corresponding families of random variables are independent for different $k$. Finally Proposition~\ref{prop:bc}(ii) together with $\sum_k\Prob(A_k)=\infty$ gives $\Prob(A_k\iio)=1$.
\end{proof}

\end{proof}
\section{Proof of Theorem~\ref{thm:rec}}\label{sec:proofrec}
\begin{proof}[Proof of Theorem~\ref{thm:rec}]
Fix $1/\beta<\alpha<\delta$ (for instance $\alpha=2/\beta$, $\delta=\alpha+1$), with $n_0$ sufficiently large. By Lemmas~\ref{lem:tail} and~\ref{lem:sausage}, the event
\[
\Omega^*:=\{D_n \text{ occurs only finitely often}\}\cap\{A_k\iio\}
\]
has probability $1$. From now on fix $\omega\in\Omega^*$ and choose a subsequence $k_1<k_2<\cdots$ with each $k_j$ large; set $K=\{k_1,k_2,\dots\}$.

\textbf{Step 1: a key inequality for the time scale.} We claim that for all sufficiently large $k$ and every $n> n_{k+1}$,
\begin{equation}\label{eq:keyineq}
n-\alpha\ln n\;>\;n_k .
\end{equation}
Indeed, the function $g(n)=n-\alpha\ln n$ is increasing for $n>\alpha$, and
\[
g(n_{k+1})=n_k+\delta\ln n_k-\alpha\ln n_{k+1}
=n_k+(\delta-\alpha)\ln n_k+\alpha\ln\frac{n_k}{n_{k+1}} .
\]
By \eqref{eq:lnnk}, when $n_k$ is sufficiently large, $\alpha\ln(n_k/n_{k+1})\geq -\alpha\delta\,\frac1{n_k}\ln n_k\geq -\tfrac12(\delta-\alpha)\ln n_k$, hence
\[
g(n_{k+1})\geq n_k+\tfrac12(\delta-\alpha)\ln n_k>n_k,
\]
and by monotonicity of $g$ we obtain \eqref{eq:keyineq}.

\textbf{Step 2: the cut consists of exactly $m$ parallel edges.} Fix $k\in K$ sufficiently large, and consider the cut $C_k$ separating $S_k:=\{v_0,\dots,v_{n_k}\}$ from its complement, i.e., the set of edges with one endpoint of index $\leq n_k$ and the other of index $>n_k$. Let $\{v_n,w_{n,i}\}$ be an edge of $C_k$. Then $n>n_k$ and $w_{n,i}$ has index $\leq n_k$. We check the three cases:
\begin{itemize}[itemsep=1mm]
\item If $n=n_k+1$: by \eqref{eq:Ak}(a), $w_{n,i}=v_{n_k}$, producing precisely the $m$ parallel edges between $v_{n_k+1}$ and $v_{n_k}$.
\item If $n_k+1<n\leq n_{k+1}$: by \eqref{eq:Ak}(b), the index of $w_{n,i}$ is $\geq n_k+1$; it produces no cut edge.
\item If $n>n_{k+1}$: by \eqref{eq:tailevent} and \eqref{eq:keyineq}, the index of $w_{n,i}$ is $>n-\alpha\ln n>n_k$; it produces no cut edge.
\end{itemize}
Therefore
\begin{equation}\label{eq:cut}
C_k=\bigl\{\,m\ \text{parallel edges connecting } v_{n_k}\ \text{and } v_{n_k+1}\,\bigr\},\qquad |C_k|=m .
\end{equation}

\textbf{Step 3: lower bound on effective resistance and recurrence.} The cuts $C_{k_1},\dots,C_{k_J}$ satisfy $S_{k_1}\subset S_{k_2}\subset\cdots$, hence any path from $v_0$ to the set $\{v_n:n>n_{k_J}\}$ must traverse each $C_{k_j}$ in turn. By the energy principle \eqref{eq:cutbound},
\[
\Reff\bigl(v_0\leftrightarrow \{v_n:n>n_{k_J}\}\bigr)\;\geq\;\sum_{j=1}^{J}\frac1{|C_{k_j}|}=\frac{J}{m}.
\]
Letting $J\to\infty$ yields $\Reff(v_0\leftrightarrow\infty)=\infty$. By Proposition~\ref{prop:electric}(i), $G(\omega)$ is recurrent. Since $\Prob(\Omega^*)=1$, the theorem is proved.
\end{proof}

\begin{remark}[Sausage-like structure]\label{rem:sausage}
The above proof actually gives a stronger structural statement: almost surely there exist infinitely many indices $k$ such that between times $n_k$ and $n_{k+1}$ the graph $G$ is ``sausage-like'' --- the entire block of new vertices only connects to the ``past'' through the $m$ parallel edges between $v_{n_k}$ and $v_{n_k+1}$. This structure is the foundation for the renewal theory of Chapter~4.
\end{remark}

\chapter{Infinite collisions}\label{ch:icp}
This chapter proves Theorem~\ref{thm:icp}: for \emph{every} $\beta>0$, $G$ has the infinite-collision property almost surely. The proof is in three steps: first, using the memoryless property of the exponential weights, we construct ``recording bottleneck times'' $\tau_k$ and the renewal process $\{T_k\}$ (\S\ref{sec:renewal}); then we show that $\tau_1$ has finite moments of every order (\S\ref{sec:moments}); finally a moderate-deviation estimate verifies the criterion in \S\ref{sec:proof-icp}.
\section{Recording bottleneck times and the renewal process}\label{sec:renewal}

From now on we fix $\beta>0$. For $n\geq 0$, define the \textbf{recording-bottleneck event}
\[
\mathcal{A}_n:=\Bigl\{\begin{array}{l}
\text{(a)}\ \forall\,1\leq i\leq m:\ w_{n+1,i}=v_n;\\[1mm]
\text{(b)}\ \forall\, j\geq n+2,\ \forall\,1\leq i\leq m:\ w_{j,i}\in\{v_{n+1},\dots,v_{j-1}\},
\end{array}\Bigr.
\]
Condition (a) says that all $m$ edges of $v_{n+1}$ are attached to $v_n$; condition (b) says that \emph{forever} afterwards no new edge returns to $\{v_0,\dots,v_n\}$. Taken together, they imply that the cut between $S_n=\{v_0,\dots,v_n\}$ and its complement consists of exactly the $m$ parallel edges between $v_n$ and $v_{n+1}$, and this cut is ``terminal'' (it will not be violated by any future edge). Set
\[
  \tau_{k+1}:=\min\{n>\tau_k:\mathcal{A}_n \text{ occurs}\},\qquad
  T_k:=\tau_k+1\ (k\ge1),\qquad T_0:=0 .
\]
$T_k$ is the right endpoint of the $k$-th ``sausage segment'': $G=\bigcup_{k\ge0}G_{T_k}$, and consecutive segments are connected only by $m$ parallel edges.

\begin{theorem}\label{thm:renewal}
$\{T_k:k\ge0\}$ is a stationary-increment renewal process on $\Z_+$, i.e.,
\[
  T_k=\sum_{j=1}^{k}X_j,\qquad X_j:=T_j-T_{j-1}\ \text{i.i.d.}\ \overset{d}{=}\ T_1=\tau_1+1,
\]
and $\Prob(T_1<\infty)=1$.
\end{theorem}

\begin{proof}
$\Prob(T_1<\infty)=1$, i.e., ``the recording time almost surely exists'', follows immediately from Lemmas~\ref{lem:D-prob} and~\ref{lem:sigma} in \S\ref{sec:moments}. The core is the \emph{i.i.d.}\ part, which relies on the memoryless property of exponential weights: we must verify, for arbitrary $j,\ell\ge1$, $1\le i\le m$ and $j\le j'\le j+\ell-1$, that
\[
  \Prob\big(w_{j+\ell,i}=v_{j'}\,\big|\,w_{j+\ell,i}\ge j\big)=\Prob\big(w_{\ell,i}=v_{j'-j}\big),
\]
where $w\ge j$ is shorthand for ``the index of $w$ is $\ge j$''. Indeed,
\[
  \frac{\Prob(w_{j+\ell,i}=v_{j'})}{\Prob(w_{j+\ell,i}\ge j)}
  =\frac{e^{\beta j'}/Z_{j+\ell}}{(e^{\beta j}+e^{\beta(j+1)}+\dots+e^{\beta(j+\ell-1)})/Z_{j+\ell}}
  =\frac{e^{\beta(j'-j)}}{Z_\ell}
  =\Prob(w_{\ell,i}=v_{j'-j}).
\]
Therefore, after ``resetting time to $T_k$'', the connection rule on $\{v_{T_k},v_{T_k+1},\dots\}$ has \emph{the same distribution} as the original process (condition (b) has already ruled out any edges returning to the past, so conditioning on ``index $\ge T_k$'' is exactly this), hence $X_{k+1}=T_{k+1}-T_k$ is independent of $(X_1,\dots,X_k)$ and identically distributed to $X_1$.
\end{proof}

\begin{remark}[Role of the memoryless property]\label{rem:memoryless}
This is the \emph{only place in the paper where exponential weights are irreplaceable}: if the connection probabilities are given by general weights $a_k$, the memoryless property generally fails and $\{T_k\}$ no longer forms a renewal process, causing the argument of Chapter~3 to collapse. The ``translation invariance'' of the probability $\big(\frac{e^{\beta j}}{\sum_{r=j}^{j+\ell-1}e^{\beta r}}\big)$ is exactly the discrete analogue of the memoryless property of the geometric distribution (or, in continuous time, of the exponential distribution).
\end{remark}
\section{Moments of $\tau_1$ and existence of a recording time}\label{sec:moments}

In this section we prove, for \emph{every} $\beta>0$, that:
\begin{proposition}\label{prop:moments}
For any $\gamma\in\N$, $\E\tau_1^{\gamma}<\infty$; in particular, $\Prob(\tau_1<\infty)=1$.
\end{proposition}

\subsection{Step 1: tail probability of ``old-vertex'' events}

Choose $\bar\alpha>1/\beta$ (in this section we may take it arbitrarily large). Set
\[
  A_n:=\Big\{\exists\,1\le i\le m:\ \text{index of } w_{n,i}\ \le n-\bar\alpha\ln n\Big\},\qquad
  B_n:=\bigcup_{\ell\ge n}A_\ell .
\]

\begin{lemma}\label{lem:tail-B}
$\Prob(A_n)\le m\,Z_{n-\bar\alpha\ln n}/Z_n\sim m\,n^{-\bar\alpha\beta}$, and there exists a constant $C=C(\bar\alpha,\beta,m)$ such that
\[
  \Prob(B_n)\ \le\ C\,n^{-(\bar\alpha\beta-1)} .
\]
In particular, by the Borel--Cantelli lemma, $\Prob(A_n\iio)=0$, i.e., almost surely the index of $w_{n,i}$ is $\ge n-\bar\alpha\ln n+1$ for all sufficiently large $n$.
\end{lemma}

\begin{proof}
The first inequality is identical to Lemma~\ref{lem:tail}. Then by a comparison with the integral,
\[
  \Prob(B_n)\le\sum_{\ell\ge n}\Prob(A_\ell)\le m\sum_{\ell\ge n}\ell^{-\bar\alpha\beta}
  \le m\int_n^\infty x^{-\bar\alpha\beta}\,dx=\frac{m}{\bar\alpha\beta-1}\,n^{-(\bar\alpha\beta-1)} .
\]
\end{proof}

\subsection{Step 2: candidate blocks and their independence}\label{sub:indep}

Choose a block-length parameter $\alpha\in(1,2)$ and $n_0$ large, and let $n_{k+1}=n_k+\lceil\alpha\ln n_k\rceil$. As in Lemma~\ref{lem:scale}(ii) we have
\begin{equation}\label{eq:nk-props}
  \frac{\ln n_k}{n_k}\le\frac1k,\qquad n_k\le M_0 k^{\alpha}\quad(k\ge1)
\end{equation}
with $M_0$ a constant. For each $k$, define the ``guaranteed time''
\[
  \bar n_k:=\min\{n:\ n-\bar\alpha\ln n\ge n_k\}\ \ge\ n_k
  \qquad\bigl(\text{easily } \bar n_k-n_k\le \bar\alpha\ln n_k+O(1)\bigr).
\]
Define the \textbf{candidate event}
\[
  D_k:=\Big\{\forall\,1\le i\le m:\ w_{n_k+1,i}=v_{n_k}\Big\}
  \ \cap\
  \Big\{\forall\, n_k+2\le n\le\bar n_k,\ \forall\,1\le i\le m:\ \text{index of } w_{n,i}\ge n_k+1\Big\}.
\]

\begin{lemma}\label{lem:D-prob}
$\Prob(D_k)\ge c(\beta)>0$ for all $k$, where
\[
  c(\beta)=(1-e^{-\beta})^{m}\prod_{s=2}^{\infty}\big(1-C_\beta e^{-\beta s}\big)^{m}>0,
\]
with $C_\beta=(1-e^{-\beta})^{-1}$ a constant.
\end{lemma}

\begin{proof}
By the independence of the variables $w_{n,i}$,
\[
  \Prob(D_k)=\Bigl(\frac{e^{\beta n_k}}{Z_{n_k+1}}\Bigr)^{m}
  \prod_{n=n_k+2}^{\bar n_k}\Bigl(1-\frac{Z_{n_k}}{Z_n}\Bigr)^{m}.
\]
The first factor satisfies
\[
  \frac{e^{\beta n_k}}{Z_{n_k+1}}=\frac{e^{\beta n_k}(e^{\beta}-1)}{e^{\beta(n_k+1)}-1}
  \ \ge\ \frac{e^{\beta n_k}(e^{\beta}-1)}{e^{\beta(n_k+1)}}=1-e^{-\beta}.
\]
For the second factor, with $n\ge n_k+2$ (write $s=n-n_k\ge 2$),
\[
  \frac{Z_{n_k}}{Z_n}=\frac{e^{\beta n_k}-1}{e^{\beta n}-1}\le\frac{e^{\beta n_k}}{e^{\beta n}-e^{\beta n}e^{-\beta n}}\le\frac{e^{-\beta s}}{1-e^{-\beta}}=C_\beta e^{-\beta s},
\]
so
\[
  \prod_{n=n_k+2}^{\bar n_k}\Bigl(1-\frac{Z_{n_k}}{Z_n}\Bigr)^{m}
  \ge\prod_{s=2}^{\bar n_k-n_k}\bigl(1-C_\beta e^{-\beta s}\bigr)^{m}
  \ge\prod_{s=2}^{\infty}\bigl(1-C_\beta e^{-\beta s}\bigr)^{m}>0,
\]
where the last inequality uses $\sum_s C_\beta e^{-\beta s}<\infty$ (an infinite product converges to a positive number iff $\sum_s C_\beta e^{-\beta s}<\infty$).
\end{proof}

Since $\bar n_k$ may exceed $n_{k+1}$ (this happens inevitably when $\bar\alpha>\alpha$), the time intervals $(n_k,\bar n_k]$ corresponding to consecutive $k$ may overlap, so $\{D_k\}$ are in general \emph{not} independent. Set
\[
  l_0:=\bigl\lceil\bar\alpha/\alpha\bigr\rceil+1 ,
\]
then for $n_0$ sufficiently large $\bar n_k\le n_{k+l_0}$ (since $\bar n_k-n_k\le\bar\alpha\ln n_k+O(1)\le l_0\alpha\ln n_k\cdot\frac{\bar\alpha}{l_0\alpha}\le n_{k+l_0}-n_k$), so the intervals $(n_{l\,l_0},\bar n_{l\,l_0}]$ ($l=0,1,2,\dots$) are pairwise disjoint, whence
\begin{equation}\label{eq:subseq-indep}
  \{D_{l\,l_0}:l\ge0\}\ \text{are mutually independent}.
\end{equation}

\subsection{Step 3: $\tau_1\le\sigma$ and tail estimates for $\sigma$}\label{sub:moment}

Define the stopping time
\[
  \sigma:=\min\big\{n_{l\,l_0}:\ D_{l\,l_0}\ \text{occurs and } B_{\bar n_{l\,l_0}}\ \text{does not occur}\big\}
  \qquad(\min\varnothing:=\infty).
\]

\begin{lemma}\label{lem:tau-le-sigma}
$\tau_1\le\sigma$ on $\{\sigma<\infty\}$.
\end{lemma}

\begin{proof}
Suppose $\sigma=n_k$ with $k=l\,l_0$. On $D_k\cap B_{\bar n_k}^{c}$ we verify that $n_k$ is a recording time:
(a)~ $w_{n_k+1,i}=v_{n_k}$, which is the first part of $D_k$;
(b)~ For $\ell\ge n_k+2$: if $\ell\le\bar n_k$, the second part of $D_k$ gives index of $w_{\ell,i}$ $\ge n_k+1$; if $\ell>\bar n_k$, by $B_{\bar n_k}^{c}$ (i.e.\ no $A_\ell$ occurs for $\ell\ge\bar n_k$),
\[
  \text{index of } w_{\ell,i}\ \ge\ \ell-\bar\alpha\ln\ell+1\ \ge\ \bar n_k-\bar\alpha\ln\bar n_k+1\ \ge\ n_k+1,
\]
where the last step uses the definition of $\bar n_k$. Hence $\tau_1\le n_k=\sigma$.
\end{proof}

\begin{lemma}\label{lem:sigma}
$\Prob(\sigma<\infty)=1$, and for every $\gamma\in\N$, $\E\sigma^{\gamma}<\infty$.
\end{lemma}

\begin{proof}
For $j\ge1$, let $k_0=k_0(j)$ satisfy $n_{k_0}<j\le n_{k_0+1}$. By \eqref{eq:nk-props}, $j\le n_{k_0+1}\le M_0(k_0+1)^{\alpha}$, hence $k_0(j)\ge (j/M_0)^{1/\alpha}-1$. If $\sigma\ge j$, then for every $l$ with $n_{l\,l_0}<j$, either $D_{l\,l_0}$ does not occur, or $B_{\bar n_{l\,l_0}}$ occurs. Restricting to indices $0\le l\le\lfloor k_0/(2l_0)\rfloor$ (for which $n_{l\,l_0}\le n_{\lfloor k_0/2\rfloor}<j$ and $\bar n_{l\,l_0}\ge\bar n_0$ is sufficiently large), we get
\[
  \Prob(\sigma\ge j)\ \le\
  \Prob\Big(B_{\bar n_{\lfloor k_0/2\rfloor}}\Big)
  +\Prob\Big(\bigcap_{0\le l\le\lfloor k_0/(2l_0)\rfloor}D_{l\,l_0}^{c}\Big)
  \ \le\ C\,\bar n_{\lfloor k_0/2\rfloor}^{-(\bar\alpha\beta-1)}
  +\bigl(1-c(\beta)\bigr)^{\lfloor k_0/(2l_0)\rfloor} .
\]
By Lemma~\ref{lem:tail-B}, $k_0(j)\ge (j/M_0)^{1/\alpha}-1$ and $\bar n_k\ge n_k\ge k$,
\[
  \Prob(\sigma\ge j)\ \le\ C'\,j^{-(\bar\alpha\beta-1)/\alpha}+C''e^{-c''j^{1/\alpha}} .
\]
Now choose $\bar\alpha>\max\bigl(1/\beta,\ (\alpha(\gamma+1)+1)/\beta\bigr)$, so that $(\bar\alpha\beta-1)/\alpha>\gamma+1$. Then
\[
  \E\sigma^{\gamma}=\gamma\int_0^\infty j^{\gamma-1}\Prob(\sigma\ge j)\,dj<\infty,
  \qquad\Prob(\sigma=\infty)\le\Prob(\sigma\ge j)\to0 .
\]
\end{proof}

\begin{proof}[Proof of Proposition~\ref{prop:moments}]
By Lemmas~\ref{lem:tau-le-sigma} and~\ref{lem:sigma}, $\tau_1\le\sigma$, $\Prob(\sigma<\infty)=1$, hence $\Prob(\tau_1<\infty)=1$ and $\E\tau_1^{\gamma}\le\E\sigma^{\gamma}<\infty$.
\end{proof}

\begin{remark}[Choice of parameters]
The order in which the parameters are chosen is: first fix $\gamma$ and $\beta$, then choose $\bar\alpha$ large enough (so that simultaneously $\bar\alpha\beta-1>\alpha(\gamma+1)$ and $\bar\alpha>1/\beta$), and finally set $l_0=\lceil\bar\alpha/\alpha\rceil+1$ to thin the block sequence and restore independence.
\end{remark}
\section{Proof of Theorem~\ref{thm:icp}}\label{sec:proof-icp}

Let $X_j=T_j-T_{j-1}$ (i.i.d., distributed as $T_1=\tau_1+1$), and $T_n=\sum_{j\le n}X_j$. By Proposition~\ref{prop:moments}, $\mu:=\E X_1<\infty$ and $\sigma^2:=\Var X_1<\infty$ (in fact, moments of every order are finite).

\begin{lemma}\label{lem:MDP}
With probability $1$, there exists $n(\omega)$ such that
\[
  T_n(\omega)\ \le\ x\,n,\qquad\forall\, n\ge n(\omega),
\]
for any fixed constant $x>\mu$.
\end{lemma}

\begin{proof}
By the moderate-deviation estimate of Proposition~\ref{prop:limit}(iii),
\[
  \Prob\bigl(T_n-n\mu\ \ge\ n^{3/4}\bigr)\ \sim\ \exp\Bigl(-\frac{n^{1/2}}{2\sigma^{2}}\Bigr),\qquad n\to\infty,
\]
and the right-hand side is summable. By the Borel--Cantelli lemma, almost surely for all sufficiently large $n$ we have $T_n\le n\mu+n^{3/4}\le xn$ (since $n^{3/4}\le(x-\mu)n$ for $n$ sufficiently large).
\end{proof}

\begin{proof}[Proof of Theorem~\ref{thm:icp}]
Fix a sample $\omega$ for which Theorem~\ref{thm:rec} and Lemma~\ref{lem:MDP} hold. Set $H_k:=G_{T_k(\omega)}$ and $\rho=v_0$. By Proposition~\ref{prop:electric}(i):
\begin{equation}
  g_{H_k}(x,x)=\Reff\bigl(x\leftrightarrow {H_k}^{c}\bigr),
\end{equation}
and the quantity that actually matters in Lemma~\ref{lem:peres} is exactly this Green function. It therefore suffices to compare effective resistances:
For $x\in H_k$, since $H_k$ is the subgraph induced by the first $T_k$ vertices, the resistance from $x$ to the complement of $H_k$ does not exceed its graph distance to the boundary,
\[
  \Reff\bigl(x\leftrightarrow H_k^{c}\bigr)\ \le\ d(x,\partial H_k)\ \le\ T_k\ \le\ x\,k ,
\]
(Lemma~\ref{lem:MDP}; $k\ge n(\omega)$; in what follows we write $x$ for this constant, avoiding confusion with the vertex, and rename it $x_0$). On the other hand, $H_k$ contains $k$ serial recording-bottlenecks, each of resistance exactly $1/m$, so
\[
  \Reff\bigl(\rho\leftrightarrow H_k^{c}\bigr)\ \ \ge\ \ \frac{k}{m},\qquad k\ge1 .
\]
Dividing:
\[
  \frac{g_{H_k}(x,x)}{g_{H_k}(\rho,\rho)}
  =\frac{\Reff(x\leftrightarrow H_k^{c})}{\Reff(\rho\leftrightarrow H_k^{c})}
  \le\ \frac{x_0\,k}{k/m}=m\,x_0,\qquad\forall\,k\ge n(\omega),\ x\in H_k .
\]
For the finitely many $k<n(\omega)$, the ratio $\sup_{x\in H_k}g_{H_k}(x,x)/g_{H_k}(\rho,\rho)$ is almost surely finite (the Green function of a finite graph is finite). Hence the constant
\[
  C(\omega)=\max\Bigl\{m\,x_0,\ \sup_{1\le k<n(\omega)}\sup_{x\in H_k}\frac{g_{H_k}(x,x)}{g_{H_k}(\rho,\rho)}\Bigr\}<\infty .
\]
By Lemma~\ref{lem:peres}, two independent simple random walks on $G(\omega)$ meet infinitely often. Since $\omega$ was arbitrary, the theorem is proved.
\end{proof}

\chapter{Summary and outlook}\label{ch:summexp}

\section{Summary}

In this paper we studied simple random walks on the exponentially age-preferential random growth graph $G(m,\beta)$, and completely proved the following results (for every $m\geq1$, $\beta>0$):

(1)~\textbf{Recurrence} (Theorem~\ref{thm:rec}). The proof reveals that $G$ is almost surely \textbf{sausage-like}: a tail estimate (Lemma~\ref{lem:tail}) together with the time-scale sequence $\{n_k\}$ (Lemma~\ref{lem:scale}) shows that the bottleneck events $A_k$ occur with uniformly positive probability along a subsequence and are independent (Lemma~\ref{lem:sausage}), so Borel--Cantelli gives infinitely many bottlenecks almost surely. Each bottleneck forms a cut of exactly $m$ parallel edges, and by Thomson's energy principle, infinitely many cuts in series make the effective resistance diverge.

(2)~\textbf{Infinite collisions} (Theorem~\ref{thm:icp}). The technical core is to exploit the memoryless property of the exponential weights: we construct a renewal process $\{T_k\}$ (Theorem~\ref{thm:renewal}), then prove that its increments are i.i.d.\ (cf.\ \S\ref{sub:indep}) and have finite moments of every order (\S\ref{sub:moment}). A moderate-deviation theorem (Proposition~\ref{prop:limit}(iii)) yields a linear upper bound on $T_k$ (Lemma~\ref{lem:MDP}), which is used to construct a random constant $C(\omega)$ verifying Peres' criterion (Lemma~\ref{lem:peres}).

The two theorems share the same geometric source --- the sausage-like structure, which both underlies recurrence and, via the renewal process, quantitatively yields the collision criterion.

\section{Outlook}

The methods and results of this paper suggest several directions for further work.

\textbf{Power-law preference.} Replace the exponential weights $e^{\beta k}$ by power-law weights $k^{-\tau}$ ($\tau>0$), i.e.,
\[
\Prob(w_{n,i}=v_k)=\frac{k^{-\tau}}{Z_n},\qquad
Z_n=\sum_{k=0}^{n-1}k^{-\tau}=O\bigl(n^{1-\tau}\bigr),
\]
Then the translational self-similarity of the model breaks down and the renewal structure is no longer directly available. However, the quantity
\[
\Reff_{\widetilde{G}_n}(v_0\leftrightarrow\partial)=\Reff\bigl(v_0,G_n^c\bigr)
\]
of the wired-boundary graph is monotone in $\tau$ (the lighter the tail, the more the weights concentrate on the newest vertices, and the larger the effective resistance), allowing one to define the critical exponent
\[
\tau_c:=\inf\Bigl\{\tau:\ \lim_{n\to\infty}\Reff\bigl(v_0,G_n^c\bigr)=\infty\Bigr\}.
\]
We expect to prove $0<\tau_c<\infty$: for $\tau>\tau_c$, $G$ is almost surely recurrent and has the infinite-collision property; for $\tau<\tau_c$, $G$ is almost surely transient. This phase transition parallels critical phenomena in long-range percolation and ``gunfighter'' graphs; a rigorous proof will likely require new methods (perhaps relying on randomly stabilised trees or continuous-time embeddings).

\textbf{Transfer of methods.} The renewal structure is essentially the probabilistic manifestation of the exponential memoryless property; it should apply to more general ``random growth processes with translational self-similarity'', such as age-preferential processes with a general rate of decay or variants of random recursive trees. The combination bottleneck--cut--energy principle, and its role in such settings, deserves a systematic investigation.

\noindent School of Mathematical Sciences, Capital Normal University, Beijing, 100048, People R China.\\
\indent Email:15213277742@163.com, wangf@cnu.edu.cn, wuxy@cnu.edu.cn
\end{document}